\documentclass[12pt,a4paper]{article}

\usepackage{amssymb,amsmath,amsthm,xspace,epsfig, amsfonts,amscd}

\usepackage[left=1.5cm, right=1.5cm, top=1.5cm]{geometry}

\begin{document}

\newcommand{\N}{\mathbb{N}}
\newcommand{\R}{\mathbb{R}}
\newcommand{\Z}{\mathbb{Z}}
\newcommand{\Q}{\mathbb{Q}}
\newcommand{\C}{\mathbb{C}}
\newcommand{\PP}{\mathbb{P}}

\newcommand{\LL}{\mathbb{L}}
\newcommand{\OO}{\mathcal{O}}

\newcommand{\esp}{\vskip .3cm \noindent}
\mathchardef\flat="115B

\newcommand{\lev}{\text{\rm Lev}}

\def\ut#1{$\underline{\text{#1}}$}
\def\CC#1{${\cal C}^{#1}$}
\def\h#1{\hat #1}
\def\t#1{\tilde #1}
\def\wt#1{\widetilde{#1}}
\def\wh#1{\widehat{#1}}
\def\wb#1{\overline{#1}}

\def\restrict#1{\bigr|_{#1}}

\def\hu#1#2{\mathsf{U}_{fin}\bigl({#1},{#2}\bigr)}
\def\ch#1#2{\left(\begin{array}{c}#1 \\ #2 \end{array}\right)}

\newtheorem{lemma}{Lemma}[section]

\newtheorem{thm}[lemma]{Theorem}
\newtheorem*{thm*}{Theorem}
\newtheorem*{lemma?}{Lemma ??}

\newtheorem{defi}[lemma]{Definition}
\newtheorem{conj}[lemma]{Conjecture}
\newtheorem{cor}[lemma]{Corollary}
\newtheorem{prop}[lemma]{Proposition}
\newtheorem{prob}[lemma]{Problem}
\newtheorem{qu}[lemma]{Question}
\newtheorem{q}[lemma]{Questions}
\newtheorem*{rem}{Remark}
\newtheorem{examples}[lemma]{Examples}
\newtheorem{example}[lemma]{Example}

\title{Collared and non-collared manifold boundaries}
\date{}
\author{Mathieu Baillif\footnote{
Haute \'ecole du paysage, d'ing\'enierie et d'architecture (HEPIA),
d\'epartement d'Informatique et syst\`emes de communications,
Gen\`eve -- Suisse. {\tt mathieu.baillif@hesge.ch}
}}
\maketitle

\abstract{\footnotesize 
We gather in this note results and examples about collared and non-collared boundaries in non-metrisable manifolds.
It is known for a long time that boundaries are collared in metrisable manifolds, but the matter
is more complicated for non-metrisable ones.
The majority of the contents of this note is well known but a bit scattered in the literature (sometimes with small inaccuracies
that we correct), and some of it is apparently not published at all
and might be new.}


\section{Introduction and definitions}

It is known since the works of M. Brown \cite[Theorems 1 \& 2]{MortonBrown:1962} in the sixties
that the boundary $\partial M$ of a topological metrisable manifold is collared in $M$, that is, there is a neighborhood $U$ of $\partial M$
in $M$ and an embedding $h:\partial M\times[0,1)\to M$ sending $\langle x,0\rangle$ to $x$ which is a homeomorphism on $U$. 
(We use brackets $\langle\,,\,\rangle$ for ordered pairs, reserving parentheses for open intervals in
ordered sets.)
Brown's result is actually more general
and shows that a locally collared closed subspace of a metrisable space is collared.
Some years later, R. Connelly  \cite{Connelly:1971} found another proof that the boundaries of metrisable manifolds are collared
with a very nice argument that works backwards in the sense that you glue the collar to the manifold and then push the manifold little by little inside the collar until
it fills it completely. The collaring embedding is then given by the inverse of this pushing.
Connelly's definition of a collar is slightly more restrictive as the embeddings are assumed to be closed (and have domain $\partial M\times[0,1]$),
and he gives details only in the case where $\partial M$ is compact, writing that the proof should work in the case $\partial M$ is strongly paracompact (see below for a definition).
Since metrisability of $M$ is a priori not necessary for Connelly's argument, with appropriate changes the proof can be adapted to the case where
$M$ is non-metrisable. D. Gauld did exactly this in his recent nice book
\cite[Theorem 3.10 \& Corollary 3.11]{Gauldbook} and showed in particular that if $\partial M$ is connected and metrisable, then it is collared in $M$, even 
if $M$ happens to be non-metrisable. 
But looking at the details, we noticed that it is {\em not} possible to guarantee that the collar
embeddings are closed (there are counter-examples),
hence the proof shows that the boundary is collared under Brown's definition but not under Connelly's. 

Despite the availability of very good texts about the general theory of non-metrisable manifolds 
(aformentioned Gauld's book \cite{Gauldbook} and the older but less elementary article by P. Nyikos in the Handbook of Set-theoretical Topology 
\cite{Nyikos:1984}, for instance), we are not aware of a reference gathering all the results and counter-examples pertaining to the collaring problem of the boundary,
hence this small note. 
It can be thought of as a convenient reference for (the few)
researchers investigating nonmetrisable manifolds and/or, for researchers in more usual areas of mathematics, as 
a catalog of the horrors awaiting those who dare to not include metrisability in the definition of 
a manifold\footnote{The author 
is a member of the first group and sees these pathologies in a similar way a fungi specialist sees the layers of mold 
which developed in a pasta dish long forgotten in the back of the fridge, but understands the urge of some mathematicians 
(the majority, actually) to flush down the toilet
the entire specimen by imposing metrisability as an unmovable feature of manifolds (and by scheduling a periodic revision of the contents of the fridge).}.
Except for small corrections to known results and some of the examples, 
the author claims no originality in the contents of this note.
As such, although we aim to remain readable even for those not used to non-metrisable manifolds,
we do not flesh out completely the details of every argument and construction since they are already available elsewhere, but 
try to convey the main ideas. 
Some familiarity with countable ordinals might be useful in Section \ref{sec:ex}.
We also tried to keep the references to a minimum,
using Gauld's and Nyikos' texts as much as possible, as well as Engelking's book \cite{Engelking} for general results in point set topology.

So, we are concerned with boundaries in manifolds, which are assumed to be Hausdorff spaces each of whose points has a neighborhood
homeomorphic to 
$\R^{n-1}\times\R_{\ge 0}$ for some (fixed) $n\ge 1$. 
If we want to emphasize the dimension, we say {\em $n$-manifold} or {\em manifold of dimension $n$}.
A {\em surface} is a $2$-manifold.
The points in the boundary $\partial M$ are those 
which do not have a neighborhood homeomorphic to $\R^n$, and $\partial M$ is itself a $n-1$-manifold without boundary and a closed subset of $M$. 
We sometimes call $M-\partial M$ the interior of $M$. The interior in the topological sense will be called ``topological interior'' to avoid ambiguity.
Manifolds share many properties with the Euclidean space, in particular they are Tychonoff, locally compact and locally connected.
Speaking of connectedness, it is usually included in the definition of manifolds when metrisability is not assumed, mainly because 
for connected manifolds a lot of properties are equivalent to metrisability (see for instance Theorem \ref{thm:manmet} below)
and the statements are more cumbersome in the general case.
But we {\em do not} assume connectedness in this note, as $\partial M$ is in general not connected, and we want it to remain a manifold.

The following definition avoids the ambiguity between Connelly's and Brown's by introducing strong collars.
As usual, $\wb{U}$ denotes the topological closure of $U$. 

\begin{defi}
   Let $B$ be a closed subset of a space $X$. We say that
   $B$ is locally collared in $X$ if there is a cover $\mathcal{U}$ of $B$ by 
   subsets of $B$ that are open in $B$ and for each $U\in\mathcal{U}$
   an embedding $h:\wb{U}\times[0,1]\to X$ which sends $\langle x,0\rangle$ to $x$ when $x\in\wb{U}$, such that 
   $h^{-1}(B)=\wb{U}\times\{0\}$ and the image of $U\times[0,1)$ is open in $X$. 
   Such an $h$ is called a local collar (of $B$ in $X$). If the embeddings are closed,
   we say that $B$ is locally strongly collared in $X$, each $h$ being a local strong collar.
   Finally, if $\mathcal{U}=\{B\}$, we say that $B$ is collared in $X$ 
   (strongly, if the embedding is closed).
\end{defi}
If the ambient space $X$ is clear, we say that $B$ is collared, omitting ``in $X$''.
Context allowing, we sometimes also call ``collar'' the image of $h$ inside $X$.
The following theorem summarises the positive results on the collaring problem.
As said above, it is essentially due to M. Brown and R. Connelly.

\begin{thm}\label{thm:main}
   Let $M$ be a connected manifold with boundary $\partial M$. Then the following hold.\\
   (a) If $\partial M$ is compact, then it is strongly collared in $M$.\\
   (b) If $\partial M$ is Lindel\"of, then it is collared in $M$.\\
   (c) If $M-\partial M$ is metrisable and $\partial M$ is collared, then $\partial M$ is Lindel\"of.\\
   (d) If $M$ is normal and $\partial M$ is countably metacompact (in particular, metrisable or Lindel\"of) and collared in $M$, then $\partial M$ is strongly collared in $M$.\\
   (e) If $M-\partial M$ is metrisable, it has a collared boundary iff it has a strongly collared boundary iff $M$ is metrisable.\\
\end{thm}

Our convention is that normal and regular spaces are Hausdorff.
The definition of a countably metacompact space is recalled in Definition \ref{def:all} below.
Recall that a {\em $0$-set} is the inverse image of $\{0\}$ for a continuous real-valued function.
In particular, $0$-sets are $G_\delta$ sets.
The connection with collars is given by the following (easy) theorem.

\begin{thm}\label{thm:submain}
   Let $M$ be a manifold with boundary $\partial M$.\\
   (a) If $\partial M$ is strongly collared in $M$, it is a $0$-set in $M$.\\
   (b) If $M-\partial M$ is metrisable, then $\partial M$ is a $0$-set in $M$.\\
   (c) If $\partial M$ is collared in $M$, it is a $G_\delta$ set.
\end{thm}

Section \ref{sec:thms} is dedicated to the proof of these theorems (and more general results from which they follow). 
Section \ref{sec:ex} contains a catalog of counter-examples (elementary, in majority) to various a priori possible generalisations.
For instance:

-- Theorem \ref{thm:main} (b) may not hold if $\partial M$ is not Lindel\"of, even if $M$ is normal 
(Example \ref{ex:octant}), even if $M-\partial M$ is metrisable and $\partial M$ connected (Example \ref{Ex:nyik}), 
even if both $M-\partial M$ and 
$\partial M$ are metrisable (Example \ref{Ex:pruf}).

-- Examples \ref{ex:demonios}--\ref{ex:demonios2} show that even if $\partial M$ is Lindel\"of and thus collared,
   it might be not strongly collared. Hence Theorem \ref{thm:main} (d) does not hold if $M$ is not normal.

-- The boundary is not a $G_\delta$ set in Example \ref{ex:octant} (which is normal) and 
   not a $0$-set in Examples \ref{ex:demonios}--\ref{ex:demonios2} (which have metrisable boundary and separable interior).
   This bad behaviour can be pushed quite far by having
   closed discrete subsets of the boundary which are not $G_\delta$ (Example \ref{ex:discr}). 

-- Also, a collared $0$-set boundary is not always strongly collared (Example \ref{Ex:c0nsc}). 

\vskip .3cm
We end this introduction by citing a nice 
result of D. Gauld and S. Greenwood \cite{GauldGreenwood:2016} (see also \cite[Proposition 3.12]{Gauldbook}).
We will not repeat the proof here, but do stress that it relies mainly on two facts:
compact components of the boundary are strongly collared (individually),
and  a metrisable manifold is a countable increasing union of connected open sets with compact closure (see Theorem \ref{thm:manmet}).
It is thus an illustration of the use of the collaring theorem for non-metrisable manifolds.
Example \ref{ex:cpctbd} shows that it is not enough to assume that $M-\partial M$ is {\em separable}, that is, has a countable dense subset.

\begin{thm}[{D. Gauld \& S. Greenwood}]\label{thm:cpctbd}
  A connected manifold $M$ with metrisable interior $M-\partial M$
  has at most countably many compact components in its boundary.
\end{thm}

This note was written after trying to answer a question on MathOverflow by Kalle Rutanen. We thank him
for the question and his remarks on a draft version of these notes, Martin W. Licht for having spotted 
a small gap in one of our proofs,
as well as David Gauld for useful email exchanges about the subject (and his kindness in general).
We also thank the referee for his/her numerous usefull remarks.


\section{Theorems}\label{sec:thms}

To avoid ambiguity, we reserve the word `boundary' only for the boundary of a manifold and use `frontier' for $\wb{U}-U$ (when $U$ is open).
We follow the set-theoretic convention of writing $\omega$ for the set of natural numbers instead of $\N$.
The following long definition contains almost every concept we need in this section.

\begin{defi}\label{def:all}
   If $B$ is a subset of the space $X$, a $B$-cover is a family of open subsets of $X$ whose union contains $B$.
   If $\mathcal{F}$ is a $B$-cover, 
   a $B$-refinement is a $B$-cover $\mathcal{G}$ such that any member of $\mathcal{G}$ 
   is contained in a member of $\mathcal{F}$. If $B=X$ we just say cover and refinement, without ``$X$-''.\\
   A subspace $B$ of $X$ is metacompact [resp. paracompact] (resp. strongly paracompact) in $X$ iff for any $B$-cover $\mathcal{F}$
   there is a $B$-refinement $\mathcal{G}$ which is point-finite  [resp. locally finite] (resp. such 
   that any member of $\mathcal{G}$ intersects only finitely many other members) as a family of subsets of $X$.
   A space is metacompact [paracompact] (strongly paracompact) if metacompact [paracompact] (strongly paracompact) in itself.\\
   A space is countably metacompact [resp. countably paracompact] 
   if any countable cover has a point-finite [resp. locally finite] refinement.\\
   A partition of unity is a family of functions $\lambda_\alpha:X\to[0,1]$, $\alpha$ in some index set $\kappa$ 
   which may be assumed to be a cardinal,
   such that for any $x\in X$ we have $\lambda_\alpha(x) = 0$ for all but finitely many $\alpha$,
   and $\sum_{\alpha\in\kappa}\lambda_\alpha(x) = 1$.
   The partition of unity is subordinate to a cover of $X$ if the support (that is, the set of points with value $\not= 0$)
   of each $\lambda_\alpha$ is included in a member of the cover.
\end{defi}
Notice that if $B$ is closed and metacompact in some space $X$, then $B$ is metacompact, and ditto for the paracompact and strongly paracompact properties. The converse does not hold as a space can be metacompact (in itself) but not metacompact in
some other space.
We can always assume that a partition of unity $\mathcal{P}$ subordinate to 
a cover $\{U_\alpha\,:\,\alpha\in\kappa\}$ 
is given by $\{\lambda_\alpha\,:\,\alpha\in\kappa\}$, where the support of $\lambda_\alpha$ lies in $U_\alpha$.
This can be arranged by first taking a choice function $c$ sending each $\lambda\in\mathcal{P}$
to some $\alpha$ such that $U_\alpha$ contains its support, and then summing those members of $\mathcal{P}$
with the same image under $c$. 
Finally, if $\alpha\in\kappa$ not in the range of $c$, define $\lambda_\alpha$ to be constant on $0$.
We shall always use such an indexing.
We will need the following classical facts due to Dieudonn\'e and Dowker.

\begin{thm}[{\cite[Remark 5.1.7, Theorems 5.1.9, 5.2.6, 5.2.8 \& Exercise 5.2.A]{Engelking}, for instance}]\label{thm:part}
   \ \\
   (a) If $X$ is a paracompact Hausdorff 
       space and $\mathcal{U}$ is a cover of $X$, then $X$ is a normal 
       space and there is a partition of unity subordinate to $\mathcal{U}$.
       \\
   (b) If~ $\mathcal{U}$ is a locally finite cover of a normal space $X$, then there is a partition of unity 
       subordinate to $\mathcal{U}$.\\
   (c) A normal space is countably paracompact iff it is countably metacompact 
       iff its product with $[0,1]$ is normal iff its product with $\R$ is normal.\\
   (d) If $\{U_\alpha\,:\,\alpha\in\kappa\}$ is a cover of a paracompact space $X$, there is 
       a locally finite refinement $\{V_\alpha\,:\,\alpha\in\kappa\}$ 
       such that $\wb{V_\alpha}\subset U_\alpha$ for each $\alpha\in\kappa$.
\end{thm}

Let us first investigate in which cases it is possible to obtain (local) strong collars from (local) collars.
First, we note that going from local collars to local strong collars is automatic in regular spaces.
This is not true for Hausdorff spaces in general (see Example \ref{ex:notlocstrong}).
\begin{lemma}\label{lemma:locstrong}
   If $B$ is closed and locally collared in the regular space $X$, then $B$ is locally strongly collared in $X$.
\end{lemma}
\proof
   Given $x\in B$ and a local collar $h:U\times[0,1)\to X$ with image the open set $W\ni x$, take an open 
   $V$ such that $x\in V\subset\wb{V}\subset W$.
   Taking a closed strip containing $\langle x,0\rangle$ inside $h^{-1}(\wb{V})$, its image is a local strong collar.
\endproof

When the collar is global, the existence of a strong collar does not automatically follow from the regularity of the space, 
something more is needed.

\begin{thm}\label{thm:normalcoll}
   Let $B$ be closed and collared in the Hausdorff space $M$. Assume that $B$ is normal (in the subspace topology) and
   countably metacompact (in itself). Then the following hold. 
   \\
   (a) If $U\subset M$ is an open neighborhood of $B$, 
       there is a collaring embedding $g:B\times[0,1]\to X$ whose image is contained in $U$.\\
   (b) If $M$ is normal, then $B$ is strongly collared.
   
\end{thm}
\proof
   $B$ is a normal countably metacompact space and is thus countably paracompact by \ref{thm:part} (c).
   Fix a collaring homeomorphism $h:B\times[0,1)\to W$, where $W$ is an open neighborhood of $B$.
   \\
   (a) 
   For each $x\in B$, fix $q_x\in(0,1)\cap\Q$ and some open $V_x\subset B$
   such that $h(V_x\times[0,q_x])\subset U$. Consider the countable cover 
   $\mathcal{G}=\{G_q \,:\,q\in\Q\}$, where $G_q = \cup_{q_x =q} V_x$,
   and let $\{\lambda_q\,:\,q\in\Q\}$ be a partition of unity subordinate to it (its existence is secured by Theorem \ref{thm:part} (b)). 
   We use the partition of unity to glue together the constant functions $x\mapsto q/2$ defined in $G_q\in\mathcal{G}$
   to obtain a continuous map, that is, we set $f(x)=\sum_{q\in\Q}\lambda_q(x) \cdot q/2$.
   Given 
   $x\in B$, 
   $\lambda_q(x)>0$ for finitely many $q$, hence $r = max\{q\in\Q\,:\,\lambda_q(x)>0\}$
   is well defined and
   $0<f(x)\le r$. It follows that
   the image by $h$ of $C=\{\langle x,t\rangle\in B\times[0,1]\,:\, t\le f(x)\}$ is contained in $U$ and
   $\langle x,t\rangle\mapsto \langle x,t/f(x)\rangle$ is a homeomorphism between $C$ and
   $B\times[0,1]$. The composition of its inverse with $h$ defines the collaring embedding $g$.
   (Notice that the image of $g$ is closed in $W$, but might not be closed in $M$.)
   \\ 
   (b) By normality of $M$ let $U$ be open such that $B\subset U\subset\wb{U}\subset W$.
   By (a), there is a collaring embedding whose image is a closed subset of $\wb{U}\subset W$,
   and hence a closed subset of $M$.
\endproof

\begin{figure}[!t]
     \begin{center}
        \epsfig{figure=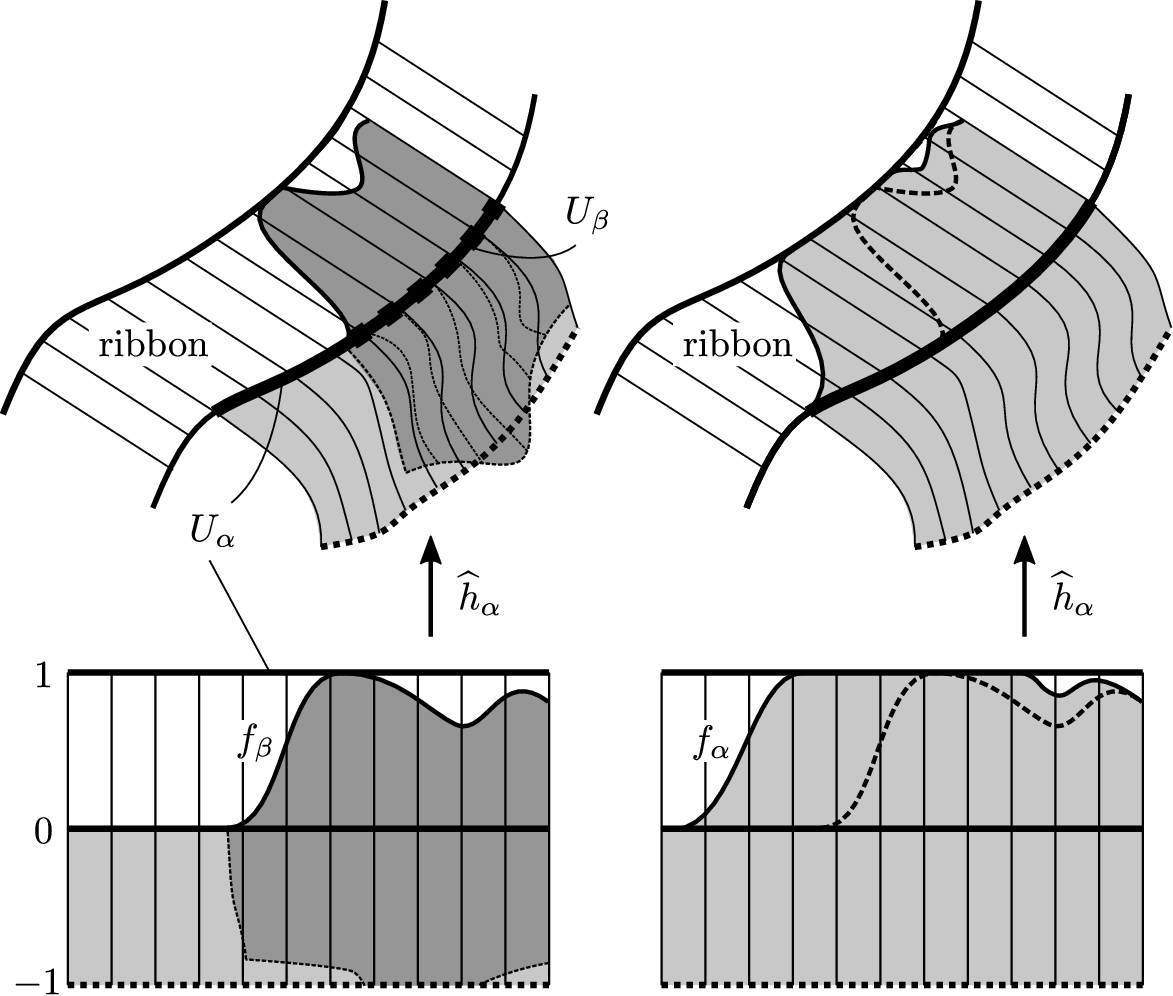, width = .90\textwidth}
     \end{center}
     \caption{Pushing in the ribbon}
     \label{fig:pushing}
   \end{figure}

Let us now turn to our main theorems.
Half of Theorem \ref{thm:main} is a direct consequence of the following more general result.
\begin{thm}[R. Connelly]
   \label{thm:Connelly}
   Let $M$ be a Hausdorff space and $B$ be a closed subset of $M$ which is locally collared in $M$ and strongly paracompact in $M$.
   Then $B$ is collared.
\end{thm}
Our proof is almost the same as Connelly's in \cite{Connelly:1971} and Gauld's version in \cite[Theorem 3.10]{Gauldbook}, 
except for some details. 
Some of them are just cosmetic and due to personal taste: we use partitions of unity, reverse some signs and try to avoid formulas,
replacing them by a geometric description.
But two details are actually somewhat important. First, we {\em do not} claim that $B$ is strongly collared, even if the local collars are strong, as
it might not be the case in general (see Examples \ref{ex:demonios}--\ref{ex:demonios2}). 
Second, we ask $B$ to be strongly paracompact in $M$ and not just in itself, as the latter
is not sufficient (see Example \ref{Ex:pruf}). 
Connelly actually only gives a complete proof when $B$ is compact and the local collars are strong
(in which case everything works fine)
and says without giving details that the proof should work in the more general case. Gauld's proof does fill the details, but 
his statement is a little bit imprecise.
\proof
   Since $B$ is locally 
   collared, we have 
   embeddings $h_\alpha:\wb{U_\alpha}\times(-1,0]\to M$ with $h_\alpha(x,0) = x$, where $\{U_\alpha\,:\,\alpha\in\kappa\}$
   is an open cover of $B$ and $\kappa$ a cardinal. Let $W_\alpha = h_\alpha(U_\alpha\times(-1,0])$.
   We first show that we can assume that
   each
   $W_\alpha$ meets only finitely many 
   other 
   $W_\beta$. \\
   By strong paracompactness of $B$ in $M$ we may take a $B$-refinement $\mathcal{W}'=\{W'_\gamma  \,:\,\gamma\in\kappa'\}$
   of $\mathcal{W}=\{W_\alpha\,:\,\alpha\in\kappa\}$ whose members intersect at most finitely other members.
   Notice that $W'_\gamma$ is a priori not the image of one of our original collaring embeddings, but we 
   now show that we can ``shrink'' our embeddings in order to have images inside these $W'_\gamma$. 
   Set $\mathcal{U}'= \{U'_\gamma \,:\,\gamma\in\kappa'\}$, where $U'_\gamma= W'_\gamma \cap B$. By Theorem
   \ref{thm:part} (d), we can shrink a bit each $U'_\gamma$ to obtain a refinement by open sets $V_\gamma$
   such that $\wb{V_\gamma}\subset U'_\gamma$.
   Notice that 
   $\wb{V_\gamma}$ is a paracompact (hence normal) space lying inside some $U_\alpha$
   and is thus collared in $M$.
   Theorem \ref{thm:normalcoll} (a) 
   provides a collaring embedding $h'_\gamma:\wb{V_\gamma}\times[0,1]\to M$ whose image lies inside $W'_\gamma$.
   Since $W'_\gamma$ intersects at most finitely many $W'_\delta$,
   the image of $h'_\gamma$ intersect the image of $h'_\delta$ for at most finitely many $\delta$.
   \\
   We may now forget about these prime symbols and 
   assume that $h_\alpha$, $W_\alpha$ and $U_\alpha$ are defined as at the beginning of the proof
   and that each
   $W_\alpha$ meets only finitely many 
   other 
   $W_\beta$.  
   We let $E=\cup_{\alpha\in\kappa}W_\alpha$.
   Since $B$ is closed and strongly paracompact in $M$, it is in particular a paracompact space, 
   hence there is a partition of unity $\{\lambda_\alpha\,:\,\alpha\in\kappa\}$
   subordinate to the cover by the $U_\alpha$'s by Theorem \ref{thm:part} (a).
   As said above we assume that the support of $\lambda_\alpha$ lies inside $U_\alpha$.
   Let $f_\alpha:B\to[0,1]$ be given by the partial sum $\sum_{\beta\le\alpha}\lambda_\beta$. Then $f_\alpha$ is continuous
   by local finiteness 
   and $f_{\kappa}(x) = 1$ for all $x\in B$.\\
   We glue the space $B\times[0,1]$ to $M$, thus defining the space $M^+$, 
   by identifying $x\in B\subset M$ to $\langle x,0\rangle\in B\times[0,1]$, and define the local embeddings 
   $\wh{h}_\alpha:B\times(-1,1]\to M^+$ which extend $h_\alpha$ 
   by the identity in $B\times[0,1]$. Call the glued part {\em the ribbon}.
   We now proceed by induction on $\alpha$ to define embeddings $\Phi_\alpha:E\to M^+$ such that
   the image of $\Phi_\alpha$ in the ribbon is $\{\langle x,t\rangle\,:\,t\le f_\alpha(x)\}$. 
   ({\em Warning:} it might be impossible to obtain an embedding from {\em all} of $M$ into $M^+$, 
   see Examples \ref{ex:demonios}--\ref{ex:demonios2} below.)
   The induction works as follows. 
   We describe geometrically what $\Phi_\alpha$ does and refer the reader more confortable with formulas 
   to the proof of \cite[Theorem 3.10]{Gauldbook} or to Connelly's original article \cite{Connelly:1971} for more details.
   We may first assume $f_{-1}$ to be $0$, $U_{-1}$ to be empty and $\Phi_{-1}$ to be the identity on $E$.
   Assume $\Phi_\beta$ is defined for each $\beta<\alpha$.
   Since $W_\alpha$ meets only finitely many other local collars 
   we may set $\beta<\alpha$ to be maximal such that $W_\alpha\cap W_\beta\not=\varnothing$
   or $-1$ if $W_\alpha$ intersects no $W_\beta$ for $\beta<\alpha$. 
   By induction the image of $\Phi_\beta$ in the ribbon is given by the points under the graph of $f_\beta$.
   Pulling back into $\wb{U_\alpha}\times[-1,1]$ with $\wh{h}_\alpha^{-1}$, we
   are in the situation depicted on the lefthandside of Figure \ref{fig:pushing}, where the image of 
   $\Phi_\beta$ inside the collar $W_\alpha$ and the ribbon is in dark grey and the rest of 
   the collar corresponding to $W_\alpha$ is depicted in lighter grey.
   To get to the righthandside of the picture, which is what we wish, it is enough to ``push along the fibers'' in
   $\wb{U_\alpha}\times(-1,1]$ using affine maps sending $(-1,f_\beta(x)]$ to $(-1,f_\alpha(x)]$ for each $x$ in $\wb{U_\alpha}$
   (and using $\wh{h}_\alpha,\wh{h}_\alpha^{-1}$ to go back and forth).
   Since $\beta<\alpha$ and we defined $f_\alpha$ and $f_\beta$ with a partition of unity, $f_\alpha$ and $f_\beta$ agree
   on the frontier of $\wb{U_\alpha}$, so this push moves only the points inside of $W_\alpha$ and defines 
   a continuous embedding of $E$ into $M^+$.
   \\
   Since each $W_\alpha$ meets only finitely many other $W_\gamma$, we can proceed with the induction until $\alpha = \kappa$,
   and local finiteness ensures the continuity of $\Phi_\alpha$ at each step.
   In the end each point of $E$ is moved only finitely many times, so letting $\Phi(y) = \Phi_\alpha(y)$
   for $\alpha$ maximal such that $y\in W_\alpha$, we obtain an embedding
   whose image contains the entire ribbon.
   The preimage of the ribbon minus $B$ is thus a homeomorphic copy of $B\times[0,1)$ in $E$.
\endproof

Recall that a space is {\em ccc} if any family of pairwise disjoint open subsets is at most countable.
A separable space is of course ccc.
\begin{thm}
   Let $X$ be a Hausdorff space and $B\subset X$ be closed.
   If $X-B$ is ccc and $B$ is not ccc, then $B$ is not collared in $X$. 
\end{thm}
\proof
   An uncountable family of disjoint open sets in $B$ yields an uncountable disjoint family of open sets in the collar.
\endproof

\begin{cor}\label{cor:comp1}
   Let $X$ be a Hausdorff ccc space and $B\subset X$ be closed and locally connected (in the subspace topology). 
   If $B$ contains uncountably many connected components, then $B$ is not collared.
\end{cor}
\proof
   Any component of $B$ is open and closed in $B$ by local connectedness (see \cite[Exercise 6.3.3]{Engelking} if in need of a proof), hence 
   $B$ is a discrete union of its components. 
\endproof

\begin{cor}\label{cor:uncollared}
   Let $M$ be a connected manifold such that $M-\partial M$ is metrisable and $\partial M$ contains uncountably many components.
   Then $\partial M$ is not collared.
\end{cor}
\proof
   A connected metrisable manifold is ccc, see just below.
\endproof

The stage is almost set for the proofs of Theorems \ref{thm:main}--\ref{thm:submain}, the only missing piece of furniture is a reminder
of the following classical properties of manifolds.
\begin{thm}[{\cite[Theorem 2.1 and A.15]{Gauldbook}}]
   \label{thm:manmet}
   Let $M$ be a connected manifold. Then the following conditions are equivalent and each imply that $M$ is separable and thus ccc. \\
   (a) $M$ is metrisable.\\
   (b) $M$ is Lindel\"of.\\
   (c) $M$ is hereditarily Lindel\"of.\\
   (d) $M$ is paracompact.\\
   (e) $M$ is strongly paracompact.\\
   (f) $M = \cup_{n\in\omega}K_n$ where $K_n$ is compact for each $n$ and the topological interior of $K_{n+1}$ contains $K_n$.\\
\end{thm}

\proof[{Proof of Theorem \ref{thm:main}}]
   \ \\
   (a) The proof of (b) below shows that $\partial M$ is collared. Since $\partial M\times[0,\frac{1}{2}]$ is compact, its homeomorphic image in $M$ is closed.\\
   (b) 
   Since $\partial M$ is Lindel\"of, we may cover it by (at most) countably many
   open sets homeomorphic to $\R^{n-1}\times\R_{\ge 0}$. The union of these open sets gives a submanifold $N\subset M$ with $\partial N =\partial M$. Since 
   $N$ is Lindel\"of it is metrisable and thus strongly paracompact by Theorem \ref{thm:manmet}. 
   Of course this implies that $\partial M$ is strongly paracompact in $N$, and we conclude with Theorem \ref{thm:Connelly}.
   \\
   (c) Since $\partial M$ is collared, there is a homeomorphic copy of $\partial M$ inside $M-\partial M$. 
       By Theorem \ref{thm:manmet}, $M-\partial M$ is hereditarily Lindel\"of, and thus so is $\partial M$.\\
   (d) This is Theorem \ref{thm:normalcoll} (b). 
       (Recall that a Lindel\"of Tychonoff space is paracompact, see e.g. \cite[Theorem 5.1.24]{Engelking}.)\\
   (e) 
       Theorem \ref{thm:manmet} tells us that Lindel\"ofness is equivalent to metrisability for connected manifolds. 
       Assume $M-\partial M$ to be metrisable. 
       Then a collared boundary implies by (c) that $M$ is metrisable, and hence normal and hereditarily paracompact,
       so the boundary is strongly collared by (d).\\
\endproof

\proof[Proof of Theorem \ref{thm:submain}]
   \
   \\
   (a) The collar is homeomorphic to $B\times[0,1]$ and closed in $M$, define 
   $f$ to be the second coordinate `in the collar' and $1$ elsewhere. 
   This defines a continuous function, because since $\partial M$ is a manifold without boundary,
   any point in the frontier of the collar is on its ``upper part'' in the interior of $M$.\\
   (b)
   Let $M$ be a manifold with metrisable interior. By Theorem \ref{thm:manmet} (f), $M-\partial M$ is
   a union $\cup_{n\in\omega}K_n$ where each $K_n$ is compact and the topological interior of $K_{n+1}$
   contains $K_n$ for each $n$.
   Hence (using either the metric or Urysohn's lemma) there is a continuous
   $f:M-\partial M\to\R_+$  such that $f(K_{n+1}-K_n)\subset[\frac{1}{n+2},\frac{1}{n+1}]$. 
   Set $\wh{f}:M\to\R_+$ to be $0$ on $\partial M$ and equal to $f$ outside of it.
   Then $\wh{f}$ is continuous and $\partial M$ is the preimage of $\{0\}$.\\
   (c) Let $h:\partial M\times[0,1)\to M$ be a collaring embedding.
       Then $\partial M$, as a subset of $M$, is the countable intersection of the open sets $h( \partial M\times[0,1/n))$.
\endproof


\section{Examples (gallery of horrors)}\label{sec:ex}

Almost all our examples are classical, although they were generally introduced in contexts 
totally unrelated to the collaring problem. We tried to give the name (if known to us)
of who first came up with the example in question (if non trivial).

\begin{example}[Trivial]
   A non-metrisable manifold with strongly collared boundary.
\end{example}
\proof[Details]
   Choose your favorite non-metrisable manifold $M$ without boundary (for example $\LL_+$ defined below). Then
   $M\times [0,1)$ trivially has a strongly collared boundary. 
\endproof

As usual, $\omega_1$ denotes the first uncountable ordinal.
Any ordinal is given by the set of its predecessors, so in particular $\omega_1$ is the set of all countable ordinals.
A countable subset of $\omega_1$ has a supremum which is a countable ordinal,
this elementary fact will be used (implicitly) many times in this section.
Recall that the longray $\LL_{\ge 0}$ is the $1$-manifold $\omega_1\times[0,1)$ endowed
with the topology given by the lexicograpic order. Its boundary contains only the point $\langle 0,0\rangle$,
removing it we obtain the open longray $\LL_+$. Chapter 1.2 in \cite{Gauldbook} is dedicated to
proving almost all the elementary properties of $\LL_+$ and $\LL_{\ge 0}$, 
recall in particular that both are normal and $\LL_{\ge 0}$ is countably compact.
We sometimes consider $\omega_1$ as a subspace of $\LL_{\ge 0}$ by identifying $\alpha\in\omega_1$ with $\langle \alpha,0\rangle\in\LL_{\ge 0}$.
The following lemma is well known.
Recall that a subset of $\omega_1$ (or $\LL_+$) is {\em stationary} if it meets any closed and unbounded (abbreviated {\em club}) 
subset of $\omega_1$  (or $\LL_+$).
A club subset of $\omega_1$ (or $\LL_+$) is stationary (see e.g. \cite[Lemma 1.15]{Gauldbook})
and contains a copy of $\omega_1$.
The definitions of stationary and club extend in an obvious way
to such sets as horizontal or vertical lines and the diagonal of $\LL_+^2$.
As one may guess, $[\alpha,\omega_1)$ is the set of points of $\LL_+$ or $\LL_{\ge 0}$ that are $\ge\alpha$,
and other intervals are defined similarly.
For any $x<y$ in $\LL_+$, the interval $(x,y)$ is homeomorphic to $\R$.

\begin{lemma}\label{lemma:propL}
   \
   \\
   (a) If $U$ is an open subset of $\LL_+^2$ whose intersection with the diagonal is stationary, 
       then
       $U$ contains $[\alpha,\omega_1)^2$ for some $\alpha\in\omega_1$.\\
   (b) If $U$ is an open subset of $\LL_+\times [0,1]$ whose intersection with the horizontal line $\LL_+\times\{a\}$
       is stationary, then $U$ contains $[\alpha,\omega_1)\times(a-\frac{1}{n},a+\frac{1}{n})$ 
       for some $\alpha\in\omega_1,n\in\omega$.
       \\
   (c) Let $C\subset\LL_+^2$ be closed.
       If $C\cap [\alpha,\omega_1)^2\not=\varnothing$ for each $\alpha\in\omega_1$, 
       then $C$ intersects the diagonal in a club
       set. If $C\cap [\alpha,\omega_1)^2=\varnothing$ for some $\alpha$, 
       then $C$ intersects either some horizontal line $\LL_+\times\{a\}$ or some vertical line $\{a\}\times\LL_+$
       in a club
       set. \\
   (d) $\LL_+\times\R$, $\LL_{\ge 0}\times\R$, $\LL_+^2$ and $\left(\LL_{\ge 0}\right)^2$ are normal spaces.
\end{lemma}
\proof
   The elementary proofs of (a) to (c) are essentially done in \cite[Lemma 3.4 \& Example 3.5]{Nyikos:1984}. 
   We now discuss the normality claims.
   It is easy to check that since $\LL_+ = (0,1)\cup[1,\omega_1)$ is the union of a Lindel\"of and a 
   countably compact space, it is countably paracompact.
   By Theorem  \ref{thm:part} (c) $\LL_+\times\R$ is normal.
   Now, assume that $A,B$ are disjoint closed sets of $\LL_+^2$.
   Recall that the intersection of two (actually, countably many) club subsets of $\LL_+$ is again club   
   \cite[Lemma 1.15]{Gauldbook}.
   The diagonal is a copy of $\LL_+$, by (c) one of $A,B$, say $B$,
   does not intersect $[\alpha,\omega_1)^2$ for some $\alpha$.
   Hence $B\subset (0,\alpha)\times\LL_+ \cup \LL_+\times(0,\alpha)$, and
   we may finish by using the result for $\LL_+\times\R$.
   The proofs for $\LL_{\ge 0}$ are almost the same.
\endproof

\begin{example}[Folklore]    
    \label{ex:octant}
    A normal surface with a boundary which is not a $G_\delta$ (and hence is non-collared).
\end{example}
\proof[Details]
  The closed octant $\mathbb{O}=\{\langle x,y\rangle\in\left(\LL_{\ge 0}\right)^2\,:\,y\le x\}$ has a boundary 
  which is not a $G_\delta$. Indeed, 
  by Lemma \ref{lemma:propL} (a) any open set containing
  the diagonal contains $\{\langle x,y\rangle\,:\,\alpha\le y\le x\}$ for some $\alpha$,
  and this remains true if one takes a countable intersection.
  As a closed subspace of a normal space, $\mathbb{O}$ is normal.
\endproof

\begin{example}[Folklore]
    A surface with a collared $0$-set boundary which is not strongly collared.
    \label{Ex:c0nsc}  
\end{example}
\proof[Details]
   Set $S = \LL_+\times[0,2) - \omega_1\times\{0\}$ (where $\omega_1$ is seen as a subset of $\LL_{\ge 0}$). 
   The boundary $(\LL_+-\omega_1)\times\{0\}$
   is an uncountable discrete union of open intervals, is collared (with collar $(\LL_+-\omega_1)\times[0,2)$) and is a $0$-set,
   but is not strongly collared. 
   To see this, let $h:(\LL_+-\omega_1)\times[0,1]\to S$ be a strong collar and $W$ be its closed image in $S$.
   Any function $\omega_1\to\R^2$ is eventually constant, meaning that it is constant in $[\alpha,\omega_1)$ for some $\alpha$,
   and a countably compact subset of a manifold is closed (see e.g. \cite[{Lemmas A.33 \& B.30}]{Gauldbook}).
   It is easily seen to imply that
   $(\LL_+-\omega_1)\times[0,2)$ does not contain a copy of $\omega_1$. Hence,
   $W$ is bounded in each horizontal line $\LL_+\times\{a\}$, $a\in(0,2)$ and 
   thus the open set $S-W$ contains $[\beta_a,\omega_1)\times(a-\frac{1}{n_a},a+\frac{1}{n_a})$ 
   for some $\beta_a,n_a$ by Lemma \ref{lemma:propL} (b).
   By Lindel\"ofness, let $E\subset\R$ be countable such that $\cup_{a\in E} (a-\frac{1}{n_a},a+\frac{1}{n_a})$ covers $(0,2)$ and
   set $\beta=\sup_{a\in E} \beta_a$. Then $S-W\supset [\beta,\omega_1)\times (0,2)$.
   This is clearly impossible as boundary components beyond $\beta$ would not be collared.
\endproof
 
\begin{figure}[!t]
     \begin{center}
        \epsfig{figure=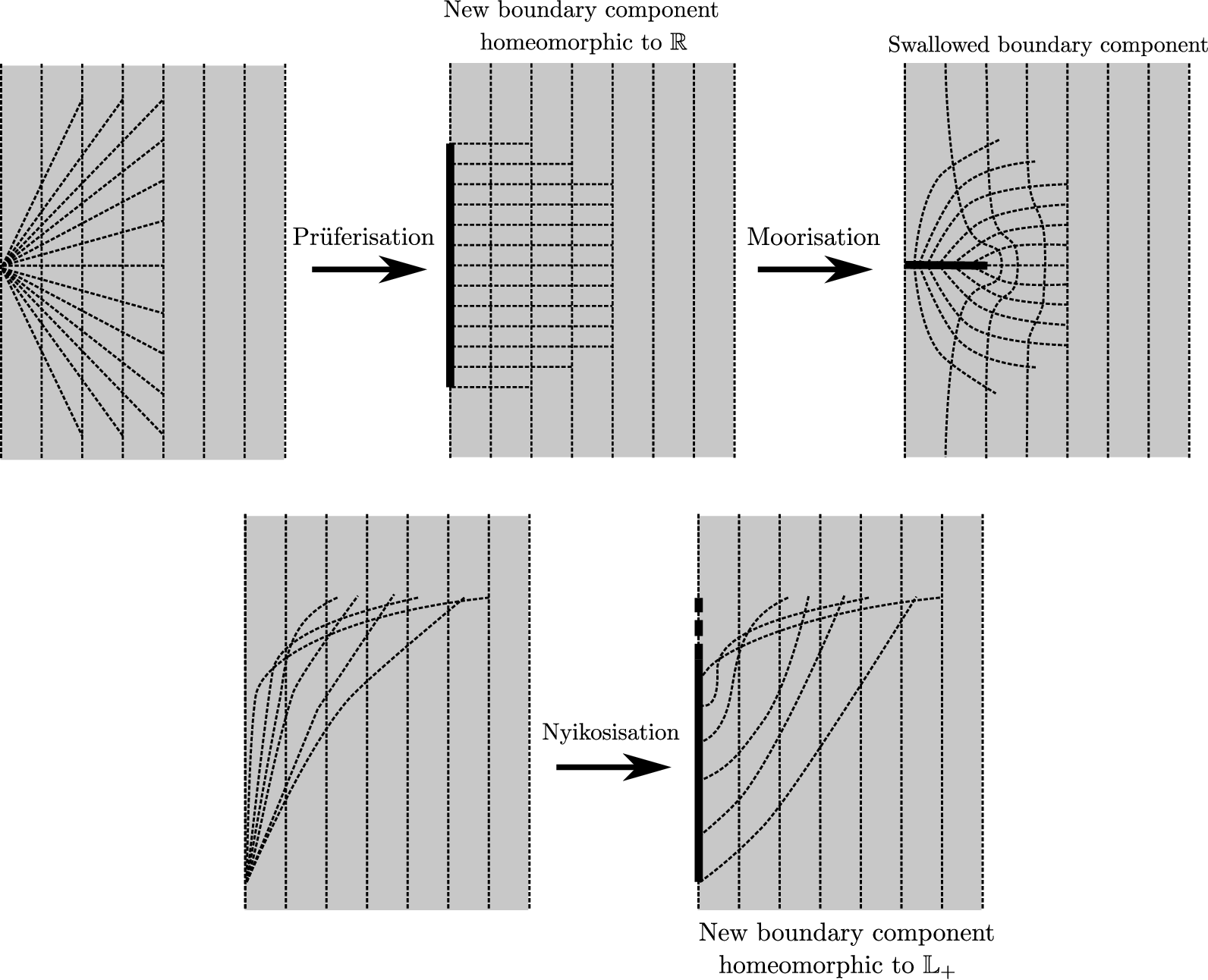, width = .9\textwidth}
     \end{center}
     \caption{Pr\"uferising and Nyikosising to glue new boundary components, Moorising to get rid of them.}
     \label{fig:prufnyik}
\end{figure}

Our next examples are based on procedures to add (or delete) boundary components in surfaces called {\em Pr\"uferisation}, {\em Moorisation} and {\em Nyikosisation}.
They are described in great detail in \cite[Chapter 1.3]{Gauldbook}, we will thus only give a geometric idea of the constructions.
\\
Start with a half plane $\R_{\ge 0}\times\R$, and choose a point $\langle 0,a\rangle$ on the vertical axis.
We replace this point by a new boundary component $\R_a$, which is a copy of $\R$, as follows.
A basic
neighborhood of some point $x\in\R_a$ is given by the union of an interval $(u,v)\subset\R_a$ containing $x$ and a triangle 
comprised between the lines of slopes $u$ and $v$ pointing at 
$\langle 0,a\rangle$ and the vertical line $\{w\}\times\R$, for some $w>0$. 
We then say that we have {\em Pr\"uferised} the half plane at height (or at point) $a$. 
Figure \ref{fig:prufnyik} is a pictural description. 
The new boundary then consists of the old boundary minus $\{\langle 0,a\rangle\}$ together with
the new boundary component $\R_a$. It takes a moment's thought to accept that we can 
perform this Pr\"uferisation at {\em each} point in the vertical axis and still have a surface
(which is the original Pr\"ufer surface, actually), because
the added boundary components are independant of each other, so to say. We can also Pr\"uferise only at some (but not all) points
and still obtain a surface, as long as 
we remove the rest of the initial boundary if 
the set of $a$ where the Pr\"ufersiation
is done is not closed.
Moreover, the added boundary components form a discrete family: 
The union of $\R_a$ and an open disk tangent to the vertical line at $\langle 0,a\rangle$ is a neighborhood
of $\R_a$ (and intersects only this boundary component).
\\
The process of {\em Moorisation} at point $a$ consists of first Pr\"uferising at $a$ and then identify $x$ with $-x$ in the added boundary component $\R_a$.
The points in $\R_a$ are in a sense pushed inside the interior of the surface.
As before, we can also Moorise at each point at once, or just some of them. \\
The {\em Nyikosisation} at point $a$ is done similarly, but this time we add a boundary component which is a copy of $\LL_+$.
Instead of taking triangles, we take a family of curves $c_\alpha$ ($\alpha\in\omega_1$) pointing at $\langle 0,a\rangle$
such that if $\beta<\alpha$, then close enough to the vertical axis $c_\alpha$ is above $c_\beta$. 
(Such a family of curves can be easily constructed from a family of functions $f_\alpha:\omega\to\omega$ such that 
$f_\beta<^* f_\alpha$ when $\beta<\alpha$, where $<^*$ means `smaller outside of a finite set'. 
See \cite[Example 1.29]{Gauldbook} for more details.)
Then a neighborhood of a point in the added longray is given by an interval in the longray union the space between 
two curves (corresponding to the endpoints of the interval) and a vertical line. 
In particular, $c_\alpha$ has $\alpha\in\LL_+$ as its unique limit point. 

Notice that the resulting surface in any of these constructions is separable, as the interior of the half plane is dense.

\begin{example}[Rado]
   A separable surface with metrisable interior and a non-collared metrisable boundary which is strongly paracompact in itself.
   \label{Ex:pruf}
\end{example}
\proof[Details]
   Pr\"uferise the half plane at each height. The resulting surface has then 
   uncountably many boundary components and thus a non-collared boundary by Corollary \ref{cor:uncollared}.
   The boundary is made of a discrete uncountable collection of copies of $\R$ and is thus metrisable
   and strongly paracompact in itself.
\endproof

\begin{example}[Nyikos]
    A separable surface with a non-collared single boundary component and metrisable interior.
    \label{Ex:nyik}  
\end{example}
\proof[Details]
   The Nyikosisation (at any point)
   of the half plane has a non-collared boundary by Theorem \ref{thm:main} (c).
\endproof

The boundary is a $0$-set in the two previous examples by Theorem \ref{thm:submain} (b).
We shall now give an example
of a manifold where it is not the case. First we need
the following fact. We were not able to spot a reference for a proof, so we give one.

\begin{lemma}\label{lemma:prufbound}
   Let $\mathbb{P}$ be the surface obtained by Pr\"uferising the positive half plane at each height
   and let $\partial \mathbb{P}=\cup_{a\in\R}\R_a$.
   Let $G\subset\R$ be a non-meagre set (hence, in particular, somewhere dense). Let $U$ be an open subset of $\mathbb{P}$ such that
   for each $a\in G$, there is some $x_a\in\R_a\cap U$. 
   Then $U$ contains a `strip' along the boundary $(0,c)\times(u,v)$ (where $c,u,v\in\R$, $c>0$, $u<v$).
\end{lemma}
\proof
   \
   \\
   \begin{minipage}{.65\textwidth}
   Let $\{q_n\,:\,n\in\omega\}$ 
   be an enumeration of $\Q$. Given $m,n,k\in\omega$ such that $q_k>q_n$ and $a\in\R$,
   let $T_{k,n,m,a}$ be the interior of the triangle in the positive half plane shown on the right.
   Given $a\in G$, there is a point $x_a$ in $U\cap\R_a$ and thus $m,n,k\in\omega$
   such that $T_{k,n,m,a}\subset U$. Let 
   $$
      G_{k,n,m} = \{a \in G\,:\, T_{k,n,m,a}\subset U\}.
   $$
   Since $G$ is the union of the $G_{k,n,m}$, one of them must be non-meagre,
   so there is an interval $(u,v)\subset\R$ in which $G_{k,n,m}$ is dense.
   Since the slopes $q_n,q_k$ and $m$ are fixed, $U$ contains a parallelogram of width $1/m$.
   The rest follows.
   \end{minipage}
   \quad
   \begin{minipage}{.3\textwidth}
     \epsfig{figure = 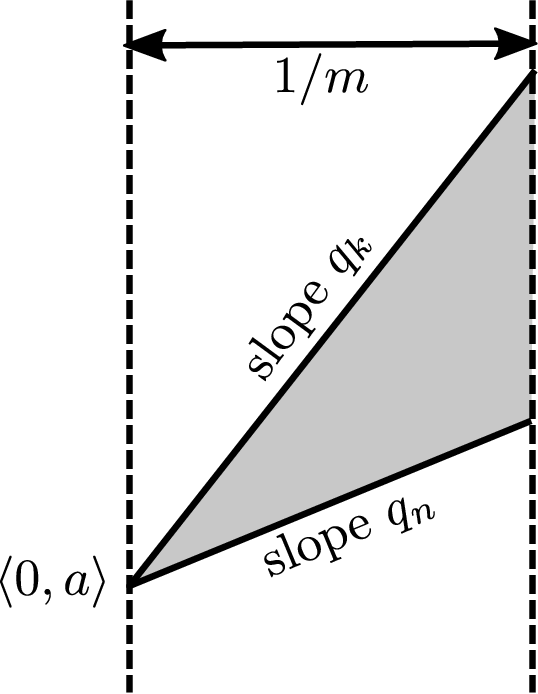, width = .6\textwidth}
   \end{minipage}
\endproof

Notice that this lemma also holds if some of the boundary components are Moorised.
\begin{example}
   A separable surface $S$ with a non $0$-set hence non-strongly collared, but metrisable and collared, boundary $\partial S$.
   Moreover, $\partial S$ is strongly paracompact in $S$.
   \label{ex:demonios}
\end{example}
\proof[Details]
   Let $S$ be the surface $S$ defined by Pr\"uferising 
   the half plane at each height and then Moorising the boundary components $\R_a$
   for $a\in\R-\Q$. The boundary being $\cup_{a\in\Q}\R_a$, it is collared by Lindel\"ofness and Theorem \ref{thm:main}.
   The argument given in the proof of Theorem \ref{thm:main} (b) above shows that Lindel\"ofness also implies that
   $\partial S$ is strongly paracompact in $S$.
   Suppose that $f:S\to[0,1]$ witnesses that $\partial S$ is a $0$-set.
   For $a\in\R-\Q$, the Moorised boundary components $\R_a$ of $\mathbb{P}$ are sent in $(0,1]$ by $f$.
   Hence there is some $n$ such that $U_n = f^{-1}((\frac{1}{n},1])$ intersects the Moorised $\R_a$ for a non-meagre set of $a$.
   By Lemma \ref{lemma:prufbound} there is some strip $(0,c)\times(u,v)$ inside of $U_n$.
   The image of $f$ on the closure of $U_n$ is contained in $[\frac{1}{n},1]$. But
   $\R_a$ is in this closure for some rational $a$, which contradicts the fact that it is a boundary component of $S$
   on which $f$ should be $0$.
\endproof

We note in passing that the (double of the) surface of Example \ref{ex:demonios} also appears in Appendix B of
\cite{FerBreVla:2021} in which D. Fern\'andez-Bret\'on and N.G. Vlamis study ends of non-metrisable manifolds.

\begin{figure}[!t]
     \begin{center}
        \epsfig{figure=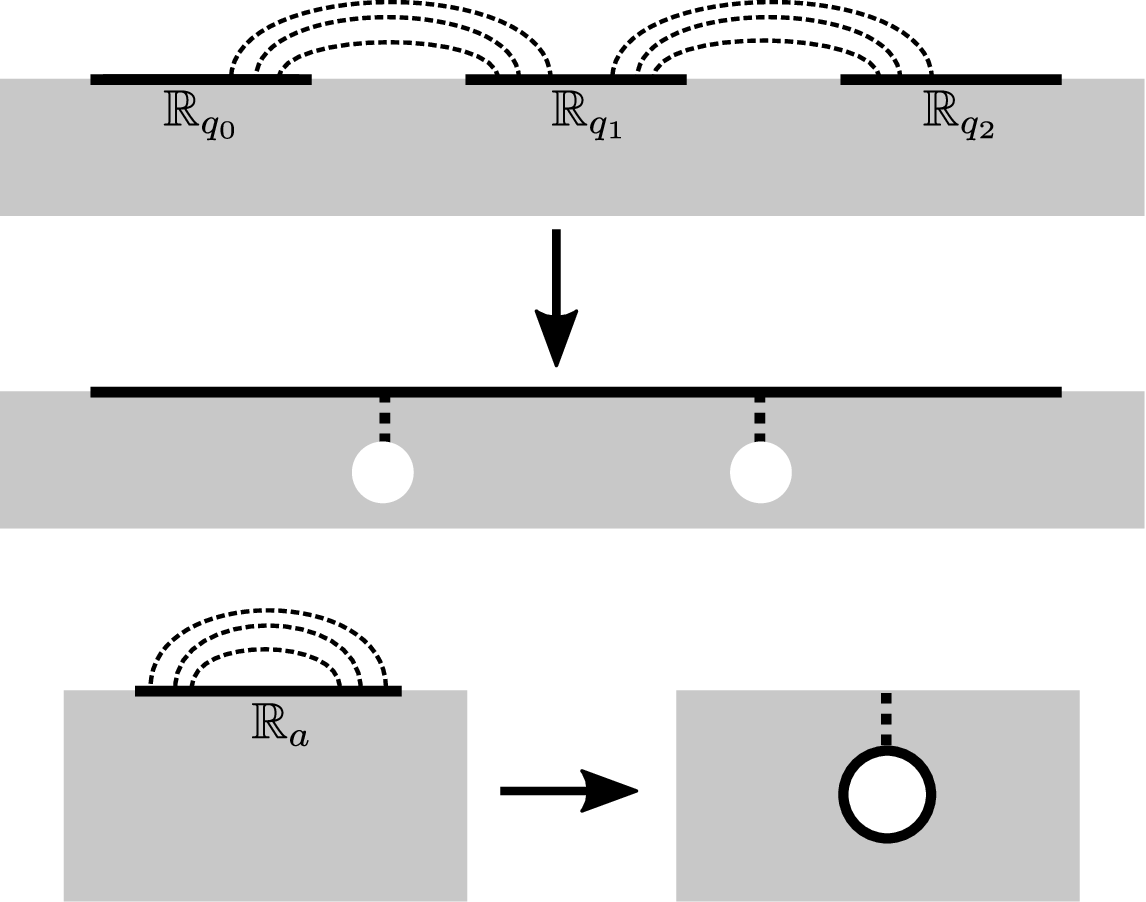, width = .80\textwidth}
     \end{center}
     \caption{Gluing boundary components in Examples \ref{ex:demonios2} (top) and \ref{ex:cpctbd} (bottom).}
        \label{fig:gluing}
   \end{figure}

\begin{example}
   A separable surface with a non $0$-set hence non-strongly collared, but metrisable, collared and connected boundary,
   which is strongly paracompact in the whole surface.
   \label{ex:demonios2}
\end{example}
\proof[Details]
   By identifying some points in the boundary in Example \ref{ex:demonios}.
   Let $q_n$, $n\in\omega$ be an enumeration of $\Q$. 
   If $B\subset\R$, we write $B_q$ for the corresponding subset of $\R_q$ (the boundary component at height $q$).
   We identify now $x\in [1,\infty)_{q_n}$ with $-x\in(-\infty,-1]_{q_{n+1}}$. 
   The resulting boundary is a copy of $\R$. See Figure \ref{fig:gluing}.
\endproof

\begin{example}[R.L. Moore, essentially]
   A separable surface with uncountably many circle boundary components.
   \label{ex:cpctbd}
\end{example}
\proof[Details]
   Take the surface of Example \ref{Ex:pruf} obtained by Pr\"uferising the half plane at each height.
   In each boundary component $\R_a$ identify $x$ with $-x$
   when $x\ge 1$. (See Figure \ref{fig:gluing}.) The non-identified points in each $\R_a$ together with $1$
   (glued to $-1$) now form
   a circle.
   The rest of $\R_a$ ends in the interior of the surface, which has thus a non-metrisable (but separable) interior.
\endproof

Let us give two last examples of manifolds with bad behaviour on the boundary by showing that even closed discrete subsets
can exhibit a reluctance to act decently.

\begin{example}\label{ex:discr0}
   The boundary of the surface of Example \ref{Ex:pruf}
   contains a countable closed discrete subset which is not a $0$-set.
\end{example}
\proof[Details]
   For each $a\in \Q$ choose one point in the boundary component $\R_a$. This defines a closed discrete set
   which is not a $0$-set, the proof being the same as the one given in Example \ref{ex:demonios}.
\endproof

However, any closed discrete subset of this surface is a $G_\delta$. 
(This can be seen roughly by first noting that any uncountable closed discrete subset is
contained in the boundary except for countably many points. Then, 
as noted in 3.5 the boundary is not collared but of course it is a $G_\delta$,
and we can argue more or less inside each boundary component. We skip the details.)

This is not the case in the next example.

\begin{example}\label{ex:discr}
   A surface with metrisable non-collared boundary which contains a closed discrete subset which is not a $G_\delta$ (hence not a $0$-set).
\end{example}
A variant of this space is cited in \cite[Remark 4.2]{KarlovaMykhailyuk:2019}, which is one of the reasons we chose to include it here.
Although we could not spot one in the literature, such examples are known for a long time,
see for instance \cite[Problem 1.2 \& Corollary 2.16]{Nyikos:1984}. (A normal example is way more elusive, since none exist
under various set-theoretic assumptions, for instance $\mathbf{V=L}$.)
The proof of the claimed properties is based on the well known Fodor's Lemma (also known as Pressing Down Lemma), 
which we recall below, and whose proof can be found in any book about set theory, for instance
\cite[II.6.15]{kunen}. 

\begin{lemma}[Fodor's Lemma]\label{lemma:fodor}
   Let $S$ be a stationary subset of $\omega_1$ and $f:S\to\omega_1$ be such that $f(\alpha)<\alpha$ for each $\alpha\in S$.
   Then there is $\alpha$ such 
   that $f^{-1}(\{\alpha\})$ is stationary.
\end{lemma}

\begin{figure}[!t]
     \begin{center}
        \epsfig{figure=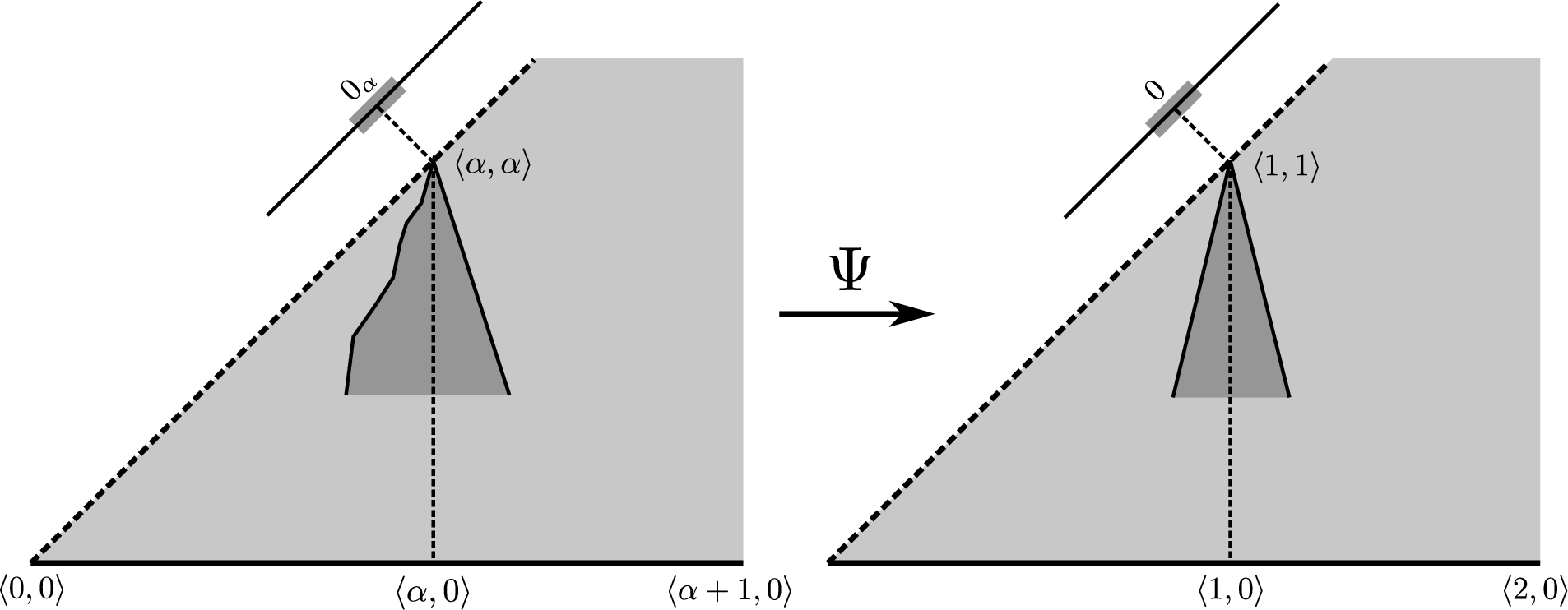, width = .9\textwidth}
     \end{center}
     \caption{Pr\"uferising on the diagonal}
     \label{fig:prufoct}
   \end{figure}

\proof[Details of Example \ref{ex:discr}]
   We start with the octant $\mathbb{O}$ of Example \ref{ex:octant}, 
   and Pr\"uferise it on the diagonal as follows. 
   The idea is to add new boundary components (real lines)
   at points $\langle\alpha,\alpha\rangle$ for $\alpha\in\omega_1$ limit. What happens to the rest of the boundary is irrelevant
   hence we might as well delete it.
   We add these components in such a way that the vertical line $\{\alpha\}\times[0,\alpha)$ 
   converges to the $0$ point of the real line attached at $\langle\alpha,\alpha\rangle$.
   We might do it as follows.
   For each limit $\alpha\in\omega_1$, choose a strictly increasing sequence $\alpha_n\in\LL_+$ converging to it.
   In what follows, the intervals $(0,\alpha),(0,\alpha]$ are understood as being in $\LL_+$ if $\alpha\ge\omega$, 
   otherwise they are 
   intervals of $\R$.
   Fix a homeomorphism $\Phi:(0,\alpha] \to (0,1]$ sending $\alpha_n$ to $1-1/n$, and extend it to a 
   homeomorphism $(0,\alpha+1) \to (0,2)$
   the obvious way.
   This yields a homeomorphism $\Psi:(0,\alpha+1)^2\to(0,2)^2$. 
   In $\{\langle x,y\rangle\in(0,2)^2\,:\,y< x\}$, Pr\"uferise at $\langle 1,1\rangle$ such that the $0$ point 
   of the added real line is the limit of the vertical line below $\langle 1,1\rangle$, and pull back this Pr\"uferisation
   in the octant by $\Psi^{-1}$. 
   Figure \ref{fig:prufoct} tries to depict what we mean by showing a 
   neighborhood of $0$ pulled back in the octant. This defines a surface, we call it $S$. 
   We write $0_\alpha$ for the $0$ point in the component added at $\langle\alpha,\alpha\rangle$.

   \begin{figure}[h]
     \begin{center}
        \epsfig{figure=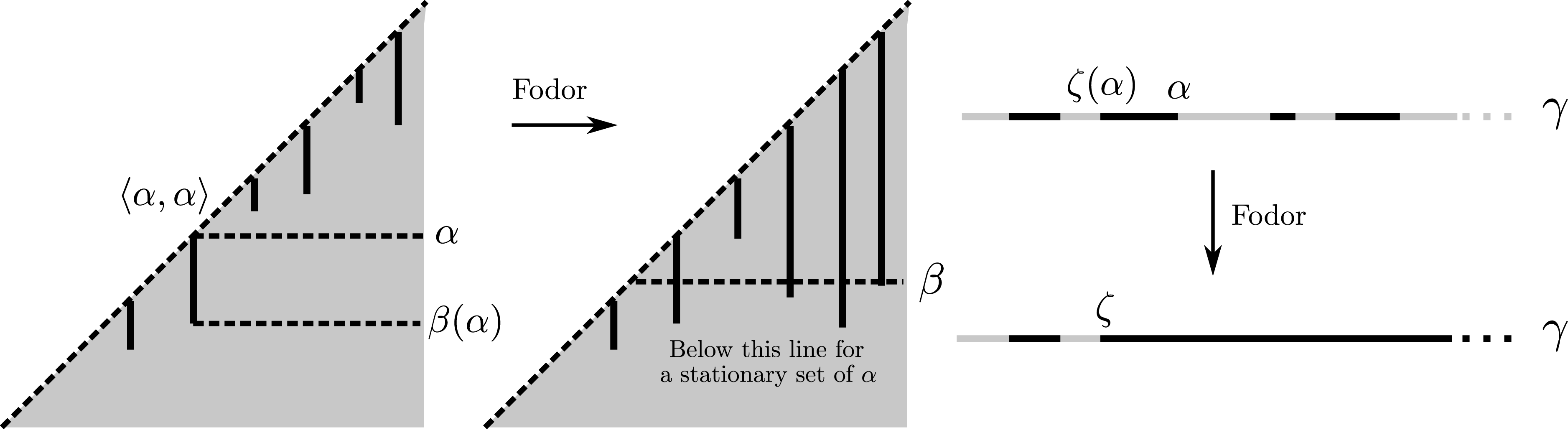, width = .9\textwidth}
     \end{center}
     \caption{Double Fodor}
     \label{fig:fodor}
   \end{figure}
   
   Then $D=\{0_\alpha\,:\,\alpha\in\omega_1,\alpha\text{ limit}\}$ is a closed discrete subset of $S$: the 
   intersection of the neighborhood depicted in Figure \ref{fig:prufoct} with $D$ is $\{0_\alpha\}$.
   We now prove that $D$ is not a $G_\delta$ by showing that for any open set $U$ containing $D$, there is $\beta\in\omega_1$
   such that for each $\gamma\ge \beta$, $U\supset [\zeta,\omega_1)\times\{\gamma\}$ for some $\zeta<\omega_1$.
   This property still holds if we take a countable intersection, proving the result.
   In passing, it also shows that the boundary is not collared.
   \\
   Let thus $U\supset D$, and for each limit $\alpha\in\omega_1$, choose $\beta(\alpha)<\alpha$ 
   such that the vertical segment $\{\alpha\}\times[\beta(\alpha),\alpha)$ lies in $U$.
   (It exists by definition of neighborhoods of $0_\alpha$.)
   By Fodor's lemma \ref{lemma:fodor}, 
   there is $\beta<\omega_1$ such that $\beta(\alpha) = \beta$ for a stationary set $E\subset\omega_1$ of $\alpha$.
   This means that if we take $\gamma\ge\beta$, $U$ 
   contains $\langle \alpha,\gamma\rangle$ for each $\alpha>\beta$ in $E$, hence for each such $\alpha$
   there is $\zeta(\alpha)<\alpha$
   with $[\zeta(\alpha),\alpha]\times\{\gamma\}\subset U$. Applying Fodor again, 
   we see that $U\supset [\zeta,\omega_1)\times\{\gamma\}$ for some $\zeta<\omega_1$, which is what we wanted to prove.\\
   Those who enjoy pictures may consult Figure \ref{fig:fodor} for a graphic 
   summary of this double application of Fodor's lemma.
\endproof

Our last example is not a manifold but shows that Lemma \ref{lemma:locstrong} does not hold when the space is not assumed to be regular.
\begin{figure}[!t]
     \begin{center}
        \epsfig{figure=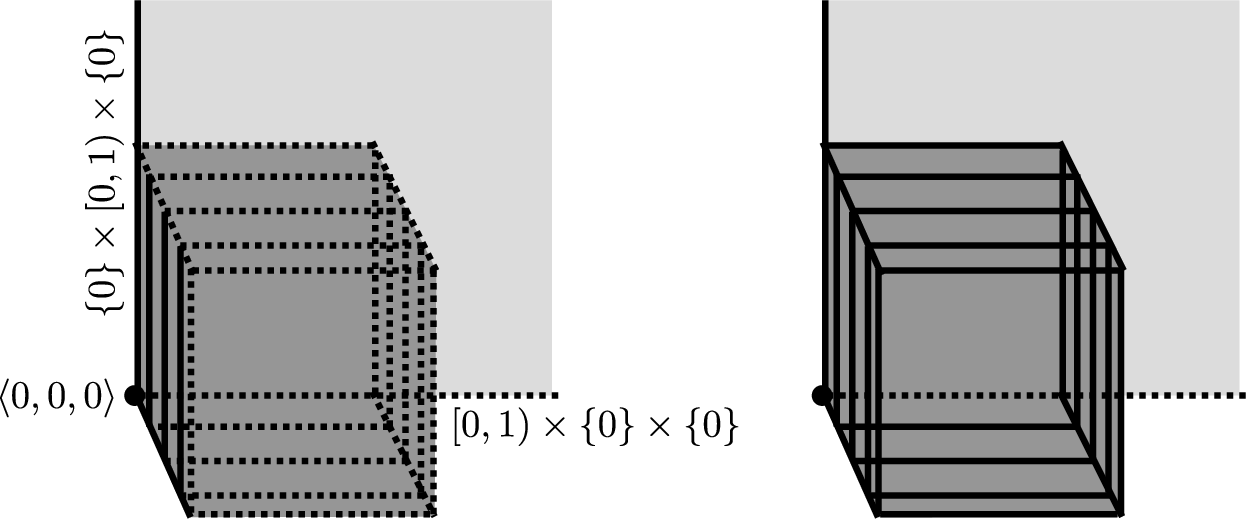, width = .7\textwidth}
     \end{center}
     \caption{A neighborhood of $\langle 0,0,0\rangle$ in Example \ref{ex:notlocstrong} (left) and its closure (right).
              The third coordinate points towards us.}
     \label{fig:nonregular}
\end{figure}

\begin{example}
   \label{ex:notlocstrong}
   A Hausdorff space $X$ with a closed locally collared subset $B$ which is not locally strongly collared.
\end{example}
\proof[Details]
   Let $H=[0,1)^2$ be given the half disk topology \cite[Example 78]{CEIT}, that is, neighborhoods of 
   points in $[0,1)\times(0,1)$ are usual open sets in the half plane
   and neighborhoods of $p=\langle x,0 \rangle$ are the union of $\{p\}$
   and an open disk centered at $p$ (or any open neighborhood of $p$ in the plane) intersected with $[0,1)\times(0,1)$.
   Then $H$ is a Hausdorff non-regular space.
   Let $X$ be $H\times[0,1)$ (endowed with product topology) with $E=(0,1)\times\{0\}\times\{0\}$ removed.
   Set $B$ to be the subset of points with $0$ third coordinate. 
   Then $B$ is closed and $\{\langle x,y,z\rangle\in X\,:\,\langle x,y,0\rangle\in B\}$ is an open collar.
   Since we removed $E$, the topology on $B$ is the usual topology as a subset of the plane, 
   in particular the product of any of its subsets with $[0,1)$ is regular.
   (Recall that any product of regular spaces is regular, see e.g. in \cite[Theorem 2.3.11]{Engelking}.)
   But any closed neighborhood of $x=\langle 0,0,0\rangle\in B$ is not regular because, as seen in Figure \ref{fig:nonregular},
   it will contain a piece of a line $[0,a]\times\{0\}\times\{c\}$ for some $a,c$.
   Hence such a closed neighborhood cannot be homeomorphic to 
   $\wb{U}\times[0,1]$ for any $U\subset B$.  There is thus 
   no closed local collar 
   at $x$. 
\endproof

\end{document}